\documentclass{article}

\usepackage[T2A]{fontenc}
\usepackage[utf8]{inputenc}
\usepackage[english]{babel}
\usepackage{amsfonts}
\usepackage{amsthm}
\usepackage{amsmath}
\usepackage{amssymb}
\usepackage{cancel}
\usepackage{tabu}
\usepackage{array}
\usepackage{tikz}
\usepackage{hyperref}
\usepackage{makecell}
\usepackage{xcolor}
\usepackage{graphicx}
\usepackage[final]{pdfpages}
\usepackage{bbm}

\usepackage[noadjust]{cite}
\usepackage{filecontents}

\graphicspath{ {./images/} }

\title{Undecidability of the elementary theory of Young--Fibonacci lattice}

\author{Vsevolod Evtushevsky}

\begin{document}

\maketitle

\begin{abstract}
For a poset $(P,\leqslant)$ we consider the first-order theory, that is defined by set $P$ and relation $\leqslant$. The problem of undecidability of combinatorial theories attracts significant attention, see for example, \cite{L1,L2,L3,L4,L5,L6,L7,L8,L9}. Recently A. Wires \cite{W} proved the undecidability of the elementary theory of Young lattice
and also established the maximal definability property
of this theory. The purpose of this article is to obtain the same results for another graded lattice, which has much in common with Young lattice: Young--Fibonacci lattice. As in \cite{W}, for the proof of undecidability we define Arithmetic into this theory.
\end{abstract}

\tableofcontents

\newpage

\newtheorem{Lemma}{Lemma}

\newtheorem{Alg}{Алгоритм}

\newtheorem{Col}{Corollary}

\newtheorem{theorem}{Theorem}

\newtheorem{Def}{Definition}

\newtheorem{Prop}{Proposition}

\newtheorem{Problem}{Задача}

\newtheorem{Zam}{Remark}

\newtheorem{Oboz}{Notation}

\newtheorem{Ex}{Пример}

\newtheorem{Nab}{Наблюдение}
\section{Introduction}

Consider the words over
the alphabet $\{1,2\}$ with
given sum of digits $n$.
It is well known that the number
of such words is Fibonacci
number $F_{n+1}$
($F_0=0,F_1=1,F_{k+2}=F_{k+1}+F_k$),
and this is the most known
combinatorial interpretation of
Fibonacci numbers.
Also such words
correspond to domino tilings
of the horizontal
rectangular strip 
$2\times n$:
digits $2$ correspond
to pairs of horizontal dominoes
and digits $1$ to vertical dominoes.

Consider a partial order
on this set of words: 
say that the word $x$ is less than or equal to
the word $y$, if, after 
removal of the maximal common suffix,
the number of digits $2$ in $y$ is not less than the total number of digits in $x$.

This relation is indeed a partial order,
and the corresponding poset is
a modular lattice known as Young --
Fibonacci lattice.

\begin{center}
\includegraphics[width=12cm, height=10cm]{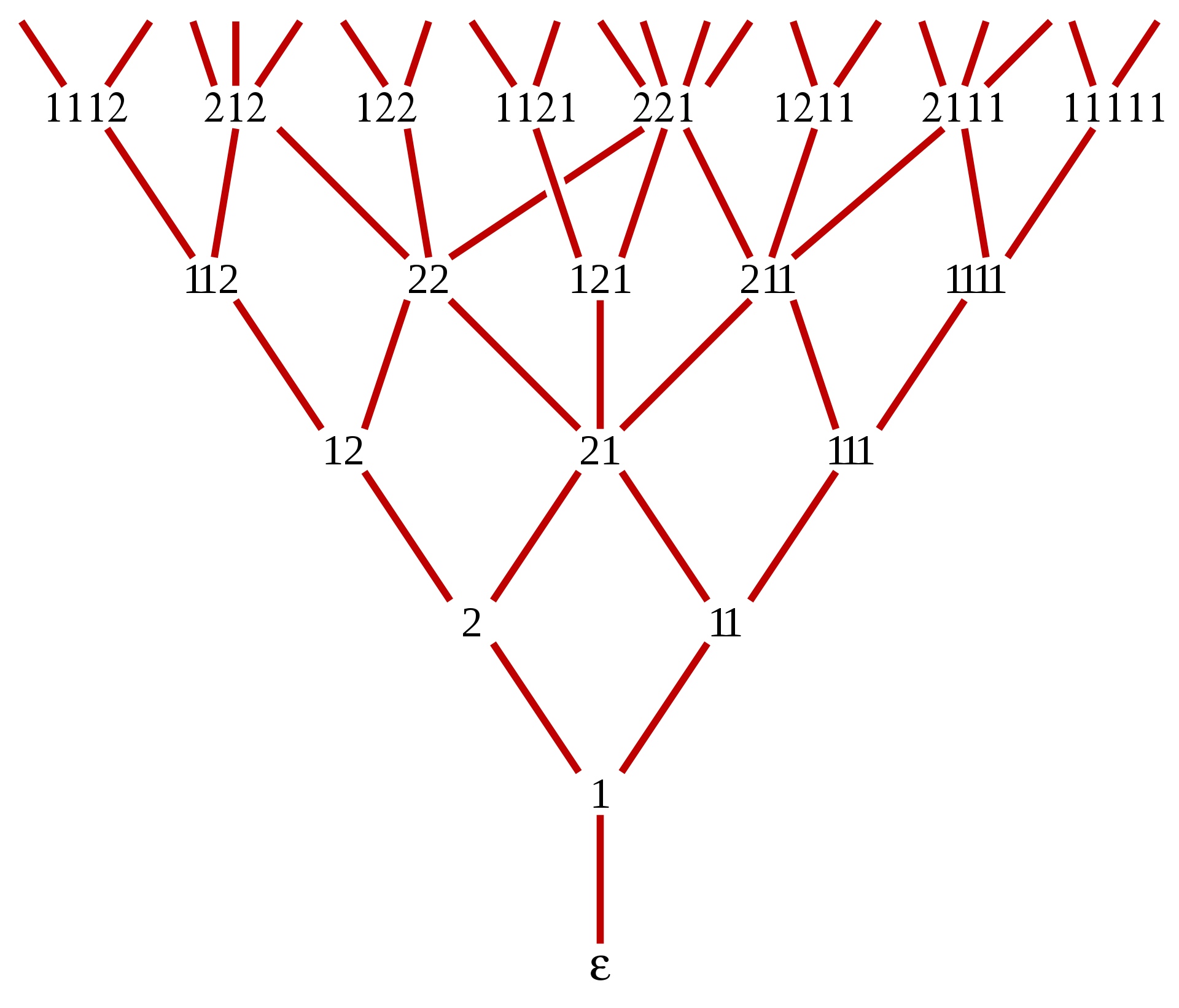}
\end{center}

Young--Fibonacci graph
(see the above figure) is 
defined as the Hasse diagram of this lattice. This is a graded
graph which grows upwards
starting from the empty word.
The grading function
is a sum of digits. The edges
go from a given word $x$ 
to the words which may be obtained from
$x$ by one of the following two
operations:
\begin{enumerate}
    \item  replace the leftmost $1$ by $2$;
    
    \item insert a $1$ anywhere to the left of the leftmost $1$.
\end{enumerate}

In addition to modularity, this graph is $1$-differential, i.e.,
for each vertex, the outdegree exceeds by $1$ the 
indegree.

The study of the graded Young--Fibonacci graph was initiated
in 1988 simultaneously and independently by such mathematicians as
Richard Stanley \cite{St} and 
Sergey Vladimirovich Fomin \cite{Fo}.

The reason why it is interesting
is that there are only two $1$-differential
modular lattices; the second one
is the lattice of Young diagrams,
which is
crucial in the theory of representations of the symmetric
group.

Recently A. Wires \cite{W} proved the undecidability of the elementary theory of Young lattice
and also established the maximal definability property
of this theory. The purpose of this article is to obtain the same results for Young--Fibonacci lattice. As in \cite{W}, for the proof of undecidability we define Arithmetic into this theory.

\newpage

\section{Notations}
\begin{Oboz}
Let $\mathbf{YF} = (\mathbb{YF},\geqslant) $ be the Young--Fibonacci lattice.
\end{Oboz}
\begin{Oboz}

Let $v \in \mathbb{YF}$. Then
    
\begin{itemize}
        \item let $|v|$ denote the sum of digits of $v$;
        \item let $\#v$ denote the number of digits of $v$;
        \item let $e(v)$ denote the number of digits $1$ in $v$;
        \item let $d(v)$ denote the number of digits $2$ in $v$.
\end{itemize}
\end{Oboz}

\begin{Zam}
Let $v\in\mathbb{YF}$. Then $|v|$ is the rank of vertex $v$ in $\mathbb{YF}$.
\end{Zam}

\begin{Oboz} Let $u,v,w\in\mathbb{YF}$. Then we write
\begin{itemize}
    \item $u=v$ instead of $u\geqslant v\land v\geqslant u$;
    \item $u\ngeqslant v$ instead of $\neg(u\geqslant v)$;
    \item $u\geqslant v \geqslant w$ instead of $u\geqslant v\land v\geqslant w$.
\end{itemize}
\end{Oboz}

\begin{Oboz}
    Let $\phi_1,\phi_2$ be first-order formulas. Then instead of
    $$ \phi_1\to\phi_2\land \phi_2\to\phi_1 $$
    we write simply
    $$\phi_1\longleftrightarrow \phi_2. $$
    
\end{Oboz}

\begin{Def}
For a structure $\left<A, \tau \right>$, a relation $R \subseteq  A^k$ is first-order definable if there is a first-order
formula $ \phi (x_1,\ldots, x_k)$ with free variables among $\{x_1,\ldots, x_k\}$ such that
$$(\pi_1,\ldots, \pi_k) \in R \Longleftrightarrow \left<A, \tau \right> \models \phi (\pi_1,\ldots,\pi_k).$$
\end{Def}

\begin{Oboz}

Let $n\in\mathbb{N}.$

    \begin{itemize}
        \item It's said, that formula belongs to set $\Pi_n,$ if it has the form 
        $$\forall y_{n,1}\ldots \forall y_{n,k_n} \exists y_{n-1,1}\ldots \exists y_{n,k_{n-1}}  \forall y_{n-2,1}\ldots \forall y_{n-2,k_{n-2}} \ldots  Q y_{1,1}\ldots Q y_{1,k_1}\phi\left(x_1,\ldots,x_k,Y\right),  $$
        where $Q=\forall$, if $n$ is odd, $Q=\exists$, if $n$ is even,  $Y$ -- are all variables that have the form $y_{i,j}$ and are quantified, and $\phi\left(x_1,\ldots,x_k,Y\right)$ is an open formula - the formula which contains no quantifiers.
        \item It's said, that the formula belongs to the set $\Sigma_n,$ if it has the form
        $$\exists y_{n,1}\ldots \exists y_{n,k_n} \forall y_{n-1,1}\ldots \forall y_{n,k_{n-1}}  \exists y_{n-2,1}\ldots \exists y_{n-2,k_{n-2}} \ldots  Q y_{1,1}\ldots Q y_{1,k_1}\phi\left(x_1,\ldots,x_k,Y\right),  $$
        where $Q=\exists$, if $n$ is odd, $Q=\forall$, if $n$ is even,  $Y$ -- are all variables that have the form $y_{i,j}$ and are quantified, and $\phi\left(x_1,\ldots,x_k,Y\right)$ is an open formula - the formula which contains no quantifiers.
    \end{itemize}
\end{Oboz}

\newpage
\section{Proof}

\begin{Prop}
    $$\left\{\varepsilon\right\}\text{-- is first-order definable in $\mathbf{YF}$}. $$
\end{Prop}
\begin{proof}
It is definable by the formula
    $$id_{\varepsilon}(u):= \forall v (v\geqslant u).$$
\end{proof}

\begin{Prop}
    $$\left\{(u,v)\in\mathbb{YF}^2:\quad \text{$u$ is a child of $v$} \right\}\text{-- is first-order definable in $\mathbf{YF}$}. $$
\end{Prop}
\begin{proof}
It is definable by the formula
    $$r(u,v):= u \geqslant v \land v \ngeqslant u\land \forall w (u\geqslant w \geqslant v \to w\geqslant u \lor v \geqslant w ).$$
\end{proof}

\begin{Oboz}
    Let $u,v\in\mathbb{YF}$, and $\phi$ be first-order formula of two variables. Then instead of
    $$ \forall w ( id_u(w) \to \phi (w,v) ) $$
    or
    $$ \exists w ( id_u(w) \land \phi (w,v) ) $$
    we write simply
    $$\phi (u,v). $$
    
\end{Oboz}

\begin{Prop}
    $$\left\{1\right\}\text{-- is first-order definable in $\mathbf{YF}$}.$$
\end{Prop}
\begin{proof}
It is definable by the formula
    $$id_1(u):=  r(u,\varepsilon). $$
\end{proof}

\begin{Prop}
    $$\left\{2,11\right\}\text{-- is first-order definable in $\mathbf{YF}$}.$$
\end{Prop}
\begin{proof}
It is definable by the formula
    $$id_{2,11}(u):=  r(u,1). $$
\end{proof}

\begin{Zam}
    There is an automorphism
    $$a: \mathbb{YF} \to \mathbb{YF},$$
    of the Young--Fibonacci lattice. It is defined as follows:
    \begin{itemize}
        \item $$\forall v\in\mathbb{YF} \quad a(v11)\to v2;$$
        \item $$\forall v\in\mathbb{YF} \quad a(v2)\to v11;$$
        \item $$\forall v\in\mathbb{YF} \quad a(v21)\to v21;$$
        \item $$a(1)=1; \quad a(\varepsilon)=\varepsilon.$$
    \end{itemize}
\end{Zam}

\begin{Col}
    In the first-order logic we can not define vertices $2$ and $11$ separately. Therefore, we add constant $2$ to our language and obtain the new structure.
\end{Col}

\begin{Oboz}
    $$\mathbf{YF}^*:=\left<\mathbb{YF},\geqslant, 2 \right>.$$
\end{Oboz}

\begin{Prop}
    $$\left\{11\right\}\text{-- is first-order definable in $\mathbf{YF}^*$}.$$
\end{Prop}
\begin{proof}
It is definable by the formula
    $$id_{11}(u):=  r(u,1)\land u\ngeqslant 2. $$
\end{proof}

\begin{Prop} $\forall u\in\mathbb{YF}$ \label{g1}
    $$\left\{u\right\}\text{-- is first-order definable in $\mathbf{YF}^*$};$$
    consequently, the automorphism $a$ is the unique nontrivial automorphism of Young--Fibonacci lattice.
\end{Prop}
\begin{proof}
For $|u|\leqslant 2$, this proposition was already proved. Further we prove proposition by induction. The base is the case $|u|\leqslant 2$. 

Induction step. Let $|u|=n\geqslant3$. Then it is easy to notice that if $v\in\mathbb{YF}:$ $|v|=n$, then sets of the parents of the vertices $u$ and $v$ can not coincide. Let $\{u_1, u_2,\ldots,u_k \}\in\mathbb{YF}^k$ be the set of parents of the vertex $u$. It is easy to notice that it is definable by the formula
    $$id_{u}(v):= \forall w \left(r(v,w)\longleftrightarrow \bigvee_{i=1}^k id_{u_i}(w)\right).$$
Now we prove that conjugation $a:\mathbb{YF} \to \mathbb{YF}$ is the unique nontrivial automorphism of  Young--Fibonacci lattice. Suppose $b:\mathbb{YF} \to \mathbb{YF}$ is another nontrivial automorphism. We know that the set $\{2, 11\}$ is first-order definable in $\mathbf{YF}$, hence it is closed under $b$.
We consider two cases:
\begin{enumerate}
    \item $b(2)=11$.
    
    In this case $b^{-1}a$ is such an automorphism that $b^{-1}a(2)=2$. It was already proved, that for all vertices $u\in\mathbb{YF}$ there exists a first-order formula $\phi_u(x,y)$, such that $u$ is the unique vertex in $\mathbb{YF}$, such that $\mathbf{YF}\models \phi_u(u,2)$. Let $R\subseteq \mathbb{YF}^2$ be the binary relation
    defined by the formula $\phi_u(x,y)$. Then $\phi_u(u,2)\in R$ implies $(b^{-1}a(u),2)=(b^{-1}a(u),b^{-1}a(2))\in R$; thus, by uniqueness of $u$ we must have $b^{-1}a(u)=u$, and so $b(u)=a(u)$. This implies $a=b$.
    \item $b(2)=2$.

    This case can be considered similarly. We should consider the automorphism $b$ instead of the automorphism $b^{-1}a$ and in the same way prove that it is trivial.
\end{enumerate}

\end{proof}
\begin{Prop}
    $$\left\{1^n:\quad n\geqslant0\right\}\text{-- is first-order definable in $\mathbf{YF}^*$} $$
\end{Prop}
\begin{proof}
It is definable by the formula
    $$\phi_{1^n}(u):=  u\ngeqslant 2. $$
\end{proof}

\begin{Prop}
    $$\left\{1^n2:\quad n\geqslant0\right\}\text{-- is first-order definable in $\mathbf{YF}^*$} $$
\begin{proof}
It is definable by the formula
    $$\phi_{1^n2}(u):= u\geqslant 2\land u\ngeqslant 11.$$
\end{proof}
\end{Prop}

\begin{Prop}
    $$\left\{1^n21:\quad n\geqslant0\right\}\text{-- is first-order definable in $\mathbf{YF}^*$} $$
\begin{proof}
It is definable by the formula
    $$\phi_{1^n21}(u):= u\geqslant 2\land u \geqslant 11\land u\ngeqslant 12 \land u\ngeqslant 111. $$
\end{proof}
\end{Prop}

\begin{Oboz}\label{o8}
Let $u,v\in\mathbb{YF}$. Then $o(u,v)\in\mathbb{YF}$ is the vertex, which we construct from $u$ and $v$ as follows:
\begin{itemize}
    \item Represent vertices $u$, $v$ in the way $u=u'w$, $v=v'w$, where $u',v',w\in\mathbb{YF}$, and $w$ is the longest common suffix of these vertices.
    \item Let $\#u\geqslant \#v$. Then in $u$ we replace the first $max(\#v'-d(u'),0)$ digits $1$ into digits $2$. 
    \item Now we get the vertex in the form of $2^xu''w$, where  $x+\#u''+\#w=\#u$, and also $x+d(u'')\geqslant\#v'$.
    \item Obviously, this notation is correct and also $o(u,v)\geqslant u$, $o(u,v)\geqslant v$.
\end{itemize}
            
\end{Oboz}
\begin{Lemma}
    Let $u,v,y\in\mathbb{YF}:$ $y\geqslant u,$ $y\geqslant v$. Then 
    $$y\geqslant o(u,v).$$
\end{Lemma}
\begin{proof}
Let $\#u\geqslant \#v$. We also introduce such designations as in the Notation \ref{o8}. Consider two cases:
\begin{enumerate}
    \item 
    Common suffix of the vertices $y$ и $u$ has at most $\#w+\#u''$ digits.

    Let $y''$ be this suffix. Then $y=y'y''$, where $y',y''\in\mathbb{YF}$. Also let $2^x u'''y''$ be $o(u,v)$, where $u'''\in\mathbb{YF}$. Then $d(y')\geqslant \#u-\#y''=x+\#u'''$, which was to be proved. 
    
    \item
    Common suffix of the vertices $y$ и $u$ has more then $\#w+\#u''$ digits.
    
    Let $y=y'u''w$, where $y'\in\mathbb{YF}$. It is easy to see that if $max(\#v'-d(u'),0)=0$, then $o(u,v)=u\geqslant v$, and in this case Lemma is obvious. And in another case, $max(\#v'-d(u'),0)=\#v'-d(u'),$ consequently $x+d(u'')=\#v'$ and also $d(y')+d(u'')\geqslant \#v'=x+d(u''),$ hence $d(y')\geqslant x$, which was to be proved. 
\end{enumerate}
\end{proof}

\begin{Prop}[corollary]
    $$\left\{(u,v,w)\in\mathbb{YF}^3: \quad w=o(u,v) \right\}\text{-- is first-order definable in $\mathbf{YF}^*$}.$$    
\end{Prop}
\begin{proof}
It is definable by the formula
    $$\phi_o(u,v,w):= w\geqslant u \land w\geqslant v \land \forall w' (w'\geqslant u \land w'\geqslant v\to w'\geqslant w). $$
\end{proof}

\begin{Zam}    
Let $n,m\in\mathbb{N}_0,$ $n\geqslant 2$. Then
\begin{itemize}
    \item if $n>m$, then
    $$o\left(1^n,1^m2\right)=2^{m+1}1^{n-m-1};$$
    \item if $n\leqslant m$, then
    $$o\left(1^n,1^m2\right)=2^{n-1}1^{m-n+1}2.$$
\end{itemize}
\end{Zam}

\begin{Oboz}
Let $n\in\mathbb{N}_0:$ $n\geqslant 2$. Then
    $$S_n:=\bigcup_{i=1}^{n-1}\left\{2^i 1^{n-i} \right\}\cup \left\{2^n\right\}\cup\bigcup_{i=1}^\infty\left\{2^{n-1}1^i2\right\}. $$

\end{Oboz}

\begin{Prop}    
$$\left\{1^n\in\mathbb{YF}: \quad n\geqslant 2\right\}\text{-- is first-order definable in $\mathbf{YF}^*$}.$$    
\end{Prop}
\begin{proof}
It is definable by the formula
    $$\phi_{1^n,11}(u):= \phi_{1^n}(u)\land 2\ngeqslant u.$$
\end{proof}

\begin{Prop}    
$$\left\{\left(1^n,v\right)\in\mathbb{YF}^2: \quad n\geqslant 2, v\in S_n \right\}\text{-- is first-order definable in $\mathbf{YF}^*$}.$$    
\end{Prop}
\begin{proof}
It is definable by the formula
    $$\phi_S(u,v):= \phi_{1^n,11}(u) \land \exists w (\phi_{1^n2}(w)\land \phi_o(u,w,v) ).$$
\end{proof}

\begin{Prop}    
$$\left\{(1^n,1^{m-2}21)\in\mathbb{YF}^2: \quad n\geqslant m \geqslant 2\right\}\text{-- is first-order definable in $\mathbf{YF}^*$}.$$    
\end{Prop}
\begin{proof}
It is definable by the formula
    $$\phi'_S(u,v):= \phi_{1^n21}(v)\land \exists w (\phi_S(u,w)\land w\geqslant v ). $$
\end{proof}

\begin{Prop}    
$$\left\{(1^n,1^{n-2}21)\in\mathbb{YF}^2: \quad n\geqslant 2\right\}\text{-- is first-order definable in $\mathbf{YF}^*$}.$$    
\end{Prop}
\begin{proof}
It is definable by the formula
    $$\phi''_S(u,v):= \phi'_S(u,v)\land \forall w (\phi'_S(u,w) \to v \geqslant w).  $$
\end{proof}

\begin{Prop}    
$$\left\{(1^n,2^{n-1}1)\in\mathbb{YF}^2: n\geqslant 2\right\}\text{-- is first-order definable in $\mathbf{YF}^*$} $$    
\end{Prop}
\begin{proof}
It is definable by the formula
    $$\phi'''_S(u,v):= \exists w (\phi''_S(u,w)\land  \phi_o(u,w,v) ).$$
\end{proof}

\begin{Prop}    
$$\left\{(1^n,2^{n})\in\mathbb{YF}^2: \quad n\geqslant 2\right\}\text{-- is first-order definable in $\mathbf{YF}^*$}.$$    
\end{Prop}
\begin{proof}
It is definable by the formula
    $$\phi'_{1^n,2^n}(u,v):=  \phi_S(u,v)\land \exists w(\phi'''_S(u,w)\land r(v,w)).$$
\end{proof}

\begin{Prop}   
$$\left\{(1^n,2^{n})\in\mathbb{YF}^2: \quad n\geqslant0\right\}\text{-- is first-order definable in $\mathbf{YF}^*$}.$$    
\end{Prop}
\begin{proof}
It is definable by the formula
    $$\phi_{1^n,2^n}(u,v):= \phi_2'(u,v) \lor id_\varepsilon(u)\land id_\varepsilon(v) \lor id_1(u)\land v=2. $$
\end{proof}

\begin{Prop}    
$$\left\{(2^n,2^{n+1})\in\mathbb{YF}^2: \quad n\in\mathbb{N}_0\right\}\text{-- is first-order definable in $\mathbf{YF}^*$}. $$    
\end{Prop}
\begin{proof}
It is definable by the formula
    $$\phi_{2^n,2^{n+1}}(u,v):= \exists w\exists w' (\phi_{1^n,2^n}(w,u) \land \phi_{1^n,2^n}(w',v) \land r(w',w) ).$$
\end{proof}

\begin{Prop}    
$$\left\{(1^n,v)\in\mathbb{YF}^2:\quad n\geqslant0, d(v)=n\right\}\text{-- is first-order definable in $\mathbf{YF}^*$}. $$    
\end{Prop}
\begin{proof}
It is definable by the formula
    $$\phi_d(u,v):= \exists w\exists w' (\phi_{1^n,2^n}(u,w)\land \phi_{2^n,2^{n+1}}(w,w') \land (v \geqslant w) \land v \ngeqslant w' ).$$
\end{proof}

\begin{Prop}    
$$\left\{(1^n,v)\in\mathbb{YF}^2:\quad n\geqslant0, \#v=n\right\}\text{-- is first-order definable in $\mathbf{YF}^*$}.$$    
\end{Prop}
\begin{proof}
It is definable by the formula
    $$\phi_{\#}(u,v):= \exists w\exists w' (\phi_{1^n,2^n}(u,w')\land \phi_{2^n,2^{n+1}}(w,w') \land w' \geqslant v \land  w \ngeqslant v ).$$
\end{proof}

\begin{Oboz}
Let $n > m\in\mathbb{N}$. Then
    $$T_{n,m}:=\bigcup_{i=1}^{m}\{2^{m-i} 1 2^i 1^{n-m-1} \}.$$
\end{Oboz}

\begin{Prop}    
$$\left\{(1^n,1^m)\in\mathbb{YF}^2: \quad n>m\geqslant 1 \right\}\text{-- is first-order definable in $\mathbf{YF}^*$}.$$    
\end{Prop}
\begin{proof}
It is definable by the formula
    $$\phi_>(u,v):= f(u)\land f(v)\land (v\ngeqslant u)\land  v\geqslant 1.$$
\end{proof}

\begin{Prop}    
$$\left\{(1^n,1^m,w)\in\mathbb{YF}^3: \quad n>m\geqslant 1, w\in T_{n,m} \right\}\text{-- is first-order definable in $\mathbf{YF}^*$}.$$    
\end{Prop}
\begin{proof}
It is definable by the formula
    $$\phi_T(u,v,w):= \phi_>(u,v) \land \phi_{\#}(u,w) \land \phi_{d}(v,w) \land  w \ngeqslant u  \land \exists u' (r (u,u') \land w\geqslant u' ).$$
\end{proof}

\begin{Oboz}
Let $n > m \in\mathbb{N}$. Then
    $$T'_{n,m}:=\bigcup_{i=1}^{m}\bigcup_{\#v'={m-i}}\{v' 1 2^i 1^{n-m-1} \}.$$
\end{Oboz}

\begin{Prop}    
$$\left\{(1^n,1^m,w)\in\mathbb{YF}^3: \quad n>m\geqslant 1, w\in T'_{n,m} \right\}\text{-- is first-order definable in $\mathbf{YF}^*$}.$$    
\end{Prop}
\begin{proof}
It is definable by the formula
    $$\phi'_T(u,v,w):= \phi_{\#}(u,w)\land \exists w' (\phi_T(u,v,w')\land (w'\geqslant w) ). $$
\end{proof}

\begin{Prop}    
$$\left\{\left(1^n,1^m,1^{m} 2 1^{n-m-1}\right)\in\mathbb{YF}^3: \quad n>m\geqslant 1 \right\}\text{-- is first-order definable in $\mathbf{YF}^*$}. $$    
\end{Prop}
\begin{proof}
    $$\phi''_T(u,v,w):= \phi'_T(u,v,w)\land \phi_d(1,w) . $$
\end{proof}

\begin{Zam}
Let $u\in\mathbb{YF}$. Then exactly one of the following conditions in true:
\begin{itemize}
    \item the leftmost digit of $u$ is $1$;
    \item the vertex $u$ has at least two parents;
    \item $u=2$;
    \item $u=\varepsilon$.
\end{itemize}
\end{Zam}

\begin{Prop}    
$$\left\{u\in\mathbb{YF}: \quad \text{the leftmost digit of $u$ is not $1$}  \right\}\text{-- is first-order definable in $\mathbf{YF}^*$}.$$    
\end{Prop}
\begin{proof}
It is definable by the formula
    $$\phi_r(u):= id_\varepsilon(u)\lor id_2(u)\lor \exists v\exists w ( r(u,v)\land r(u,w) \land  u\ngeqslant v ).$$
\end{proof}

\begin{Prop}    
$$\left\{(u,v)\in\mathbb{YF}^2: \quad \text{$v$ can be constructed from $u$ by removing some prefix of digits $1$}  \right\}$$
$$\text{-- is first-order definable in $\mathbf{YF}^*$}.$$    
\end{Prop}
\begin{proof}
It is definable by the formula
    $$\phi'_e(u,v):= u\geqslant v \land \exists w\exists w' (\phi_d(w,u)\land \phi_d(w',v)\land w = w') .$$
\end{proof}

\begin{Prop}    
$$\left\{(u,v)\in\mathbb{YF}^2: \quad \text{$v$ can be constructed from $u$ by removing the longest prefix of digits $1$}  \right\}$$
$$\text{-- is first-order definable in $\mathbf{YF}^*$}.$$    
\end{Prop}
\begin{proof}
It is definable by the formula
    $$\phi_e(u,v):= \phi_r(v)\land \phi'_e(u,v)\land \forall w (\phi_r(w)\land \phi'_e(u,w)\to v\geqslant w).$$
\end{proof}

\begin{Prop}    
$$\left\{\left(1^n,1^m, 2 1^{n-m-1}\right)\in\mathbb{YF}^3: \quad n>m\geqslant 1 \right\}\text{-- is first-order definable in $\mathbf{YF}^*$}.$$    
\end{Prop}
\begin{proof}
It is definable by the formula
    $$\phi'''_T(u,v,w):=  \exists w' (\phi''_T(u,v,w')\land \phi_e (w',w)). $$
\end{proof}

\begin{Prop}    
$$\left\{\left(1^n,1^m, 1^l\right)\in\mathbb{YF}^3:\quad n>m\geqslant 1, l\geqslant0, n=m+l \right\}\text{-- is first-order definable in $\mathbf{YF}^*$}.$$    
\end{Prop}
\begin{proof}
It is definable by the formula
    $$\phi'_+(u,v,w):= \exists  w'(\phi'''_T(u,v,w')\land  \phi_\#(w,w')). $$
\end{proof}

\begin{Prop}    \label{g2}
$$\left\{(1^n,1^m, 1^l)\in\mathbb{YF}^3:\quad n,m,l\geqslant0, n=m+l \right\}\text{-- is first-order definable in $\mathbf{YF}^*$}.$$    
\end{Prop}
\begin{proof}
It is definable by the formula
    $$\phi_+(u,v,w):= \phi'_+(u,v,w) \lor (\phi_{1^n}(u)\land  u = v \land id_\varepsilon (w)) \lor (\phi_{1^n}(u)\land u = w \land id_\varepsilon (v)) .$$
\end{proof}

\begin{Prop}    
$$\left\{(1^n,v)\in\mathbb{YF}^2:\quad n\geqslant 0, v=1^n2v' \text{ или } v=1^n \right\}\text{-- is first-order definable in $\mathbf{YF}^*$}. $$    
\end{Prop}
\begin{proof}
It is definable by the formula
    $$\phi_R(u,v):= \exists w\exists w' \exists v'( \phi_e(v,v')\land \phi_\#(w,v)\land  \phi_\#(w',v')\land \phi_+(w,w',u)  ).$$
\end{proof}

\begin{Oboz}
    Let $n\geqslant m\geqslant 0$. Then let the set of vertices that contain at most $m$ digits $2$ and whose predecessors have length of prefixes of digits $1$ at most $n$ be denoted as
    $$R_{n,m}.$$
\end{Oboz}

\begin{Oboz}
$$\phi_{\geqslant}(u,v):=\phi_{1^n}(u)\land \phi_{1^n}(v) \land u\geqslant v.$$
\end{Oboz}

\begin{Prop}    
$$\left\{(1^n,1^m,w)\in\mathbb{YF}^3:\quad n\geqslant m \geqslant 0, w\in R_{n,m}  \right\}\text{-- is first-order definable in $\mathbf{YF}^*$}.$$    
\end{Prop}
\begin{proof}
It is definable by the formula
    $$\phi'_R(u,v,w):=\phi_\geqslant (u,v)\land \exists v' (\phi_d(v',w)\land v\geqslant v')\land \forall u'\forall w' (w\geqslant w'\land \phi_R(u',w')\to u\geqslant u') . $$
\end{proof}

\begin{Zam}
    Let $n\geqslant m\geqslant 0$. Then the longest vertex that contains at most $m$ digits $2$ and whose predecessors have length of prefixes of digits $1$ at most $n$ has the form of
    $$1^n21^{n-1}2\ldots 1^{n-m+1}2 1^{n-m}.$$
\end{Zam}

\begin{Prop}    
$$\left\{\left(1^n,1^m,1^n21^{n-1}2\ldots 1^{n-m+1}2 1^{n-m}\right)\in\mathbb{YF}^3:\quad n\geqslant m \geqslant 0  \right\}\text{-- is first-order definable in $\mathbf{YF}^*$}.$$    
\end{Prop}
\begin{proof}
It is definable by the formula
    $$\phi''_R(u,v,w):=\phi'_R(u,v,w)\land \forall w'\forall w''\forall w'''(\phi'_P(u,v,w')\land \phi_\#(w'',w)\land \phi_\#(w''',w') \to w''\geqslant w''' ).$$
\end{proof}

\begin{Zam}
    $$\#(1^n21^{n-1}2\ldots 1^{n-m+1}2 1^{n-m})=n+mn-\frac{m^2-m}{2}.$$

\end{Zam}

\begin{Prop}    
$$\left\{\left(1^n,1^m,1^{nm-1/2(m^2-m)}\right)\in\mathbb{YF}^3:\quad n\geqslant m \geqslant 0  \right\}\text{-- is first-order definable in $\mathbf{YF}^*$}.$$    
\end{Prop}
\begin{proof}
It is definable by the formula
    $$\phi'_\times(u,v,w):=  \exists w'\exists w'' (\phi''_R(u,v,w')\land \phi_\#(w'',w') \land \phi_+(w'',u,w)).$$
\end{proof}

\begin{Prop}    
$$\left\{\left(1^m,1^{1/2(m^2-m)}\right)\in\mathbb{YF}^2:\quad m \geqslant0  \right\}\text{-- is first-order definable in $\mathbf{YF}^*$}.$$    
\end{Prop}
\begin{proof}
It is definable by the formula
    $$\phi''_\times(u,v):=   \exists w 
    ( \phi'_\times(u,u,w) \land \phi_+(w,u,v)) .$$
\end{proof}

\begin{Prop}    
$$\left\{(1^n,1^m,1^{l})\in\mathbb{YF}^3:\quad n\geqslant m\geqslant0, l\geqslant 0, nm=l \right\}\text{-- is first-order definable in $\mathbf{YF}^*$}. $$    
\end{Prop}
\begin{proof}
It is definable by the formula
    $$\phi'''_\times(u,v,w):=  \exists w'\exists w'' (\phi'_\times(u,v,w')\land \phi''_\times(v,w'')\land \phi_+(w,w',w'') ). $$
\end{proof}

\begin{Prop}    \label{g3}
$$\left\{\left(1^n,1^m,1^{l}\right)\in\mathbb{YF}^3:\quad n,m,l\geqslant0, nm=l  \right\}\text{-- is first-order definable in $\mathbf{YF}^*$}.$$    
\end{Prop}
\begin{proof}
It is definable by the formula
    $$\phi_\times(u,v,w):= \phi'''_\times(u,v,w)\lor \phi'''_\times(v,u,w).$$
\end{proof}

\begin{Zam}
    It follows from the proof that there exists such $m\in\mathbb{N}$, that $\phi_+(u,v,w)\in\Pi_m$ and $\phi_\times(u,v,w)\in\Pi_m.$ 
\end{Zam}

Let $\left<\mathbb{N}_0,+,\times\right>$ denote the structure over the set of non-negative integers such that the operations of addition and multiplication have their usual meaning. Propositions \ref{g1}, \ref{g2} and \ref{g3} state that we have an interpretation of $\left<\mathbb{N}_0,+,\times\right>$ into $\mathbf{YF}^*$ in which the ternary relations for addition and multiplication are definable over the vertices of form $1^n$ by $\Pi_m$-formulas. Undecidability of the positive $\Sigma_1$-theory of $\left<\mathbb{N}_0,+,\times\right>$ established in Matiyasevich \cite{Ma} yields the following:

\begin{theorem}
    The $\Sigma_{m+1}$-theory of $\mathbf{YF}^*$ is undecidable.
\end{theorem}

Since the elementary theory of a fixed structure is complete, by \cite{Tar} (Theorem 1, Theorem 7, Theorem 10) the above interpretation establishes the following:

\begin{theorem}
The elementary theory of Young--Fibonacci lattice is undecidable and non-finitely axiomatizable.
\end{theorem}

Also we want to prove that structure $\left<\mathbb{YF},\geqslant, 2 \right>$ has the maximal definability property. From the paper of Alexander Wires\cite{W}(section 4) it follows that the proof comes down to the building of a bijection between sets $\mathbb{YF}$ and $\mathbb{N}_0$ (that is a bijection $b$ from Notation \ref{bi}) and proving Proposition \ref{final} concerning this bijection.

\begin{Oboz}
    Let $\{p_0,p_1,p_2,\ldots\}=\{2,3,5,\ldots\}$ be the sequence of primes.
\end{Oboz}

\begin{Oboz}\label{bi}
    We introduce a bijection $b:\mathbb{YF}\to \mathbb{N}_0$ as follows:
    \begin{equation*}
b(v)= 
 \begin{cases}
   $$0$$ &\text {if $v = \varepsilon$}\\
   2^{n-1} &\text{if $v=1^n,\; n\geqslant1$}\\
   p_{d(v)}^{e_{d(v)}+1}\cdot\prod_{i=1}^{d(v)-1} p_i^{e_i} &\text{if $d(v)\geqslant1$, $v=1^{e_0}21^{e_1}2\ldots 2 1^{e_{d(v)-1}}21^{e_{d(v)}} .$}
 \end{cases}
\end{equation*}

\end{Oboz}

Also it follows from the paper\cite{W}(section 4), that
\begin{Prop}\label{primexp}
$$\left\{(1^n,1^m,1^l)\in\mathbb{YF}^3:  \quad n\geqslant 0, m,l\geqslant 1, \text{$p_n$ appears in the prime factorization of $l$ with exponent $m$}  \right\}$$
$$\text{-- is first-order definable in $\mathbf{YF}^*$}. $$        
\end{Prop}

\begin{Oboz}
    Let the formula that defines the set in Proposition \ref{primexp}  be denoted as 
    $$\psi_{Primexp}(u,v,w).$$
\end{Oboz}

\begin{Prop}    
$$\left\{(1^n,v,w)\in\mathbb{YF}^3:  \quad n\geqslant 0, w\geqslant v, d(w)=d(v)+n  \right\}\text{-- is first-order definable in $\mathbf{YF}^*$}. $$    
\end{Prop}
\begin{proof}
It is definable by the formula
    $$\phi_{d+}(u,v,w):=w\geqslant v\land \exists v'\exists w'( \phi_d(v',v) \land \phi_d(w',w)\land \phi_+(w',v',u)) .$$
\end{proof}
\begin{Prop}    
$$\left\{(1^n,v,w)\in\mathbb{YF}^3: \quad  n\geqslant 0, w\geqslant v, d(w)=d(v)+n, w=1^{e_0}21^{e_1}2\ldots 2 1^{e_{d(w)-1}}21^{e_{d(w)}}, e(v)=n+\sum_{i=n}^{d(w)}e_i  \right\}$$
$$\text{-- is first-order definable in $\mathbf{YF}^*$}. $$    
\end{Prop}
\begin{proof}
It is definable by the formula
    $$\phi_{d+}'(u,v,w):= \phi_{d+}(u,v,w) \land \forall v'\forall v''\forall v'''(\phi_{d+}(u,v',w)\land \phi_\#(v'',v) \land \phi_\#(v''',v')\to v''\geqslant v''' ) .$$
\end{proof}

\begin{Prop}    
$$\left\{(1^n,v)\in\mathbb{YF}^2: \quad n\geqslant0, \text{length of maximal prefix of digits $1$ of $v$ is exactly $n$}  \right\}$$
$$\text{-- is first-order definable in $\mathbf{YF}^*$}.$$    
\end{Prop}
\begin{proof}
It is definable by the formula
    $$\phi_E(u,v):= \forall w\forall v'\forall w' (\phi_e(v,w)\land \phi_\#(v',v) \land \phi_\#(w',w) \to \phi_+(v',w',u)).$$
\end{proof}

\begin{Prop}    
$$\left\{(1^n,v,w)\in\mathbb{YF}^3: \quad d(w) \geqslant n\geqslant 0, w=1^{e_0}21^{e_1}2\ldots 2 1^{e_{d(w)-1}}21^{e_{d(w)}}, v=1^{n+e_n}21^{e_{n+1}}2\ldots 2 1^{e_{d(w)-1}}21^{e_{d(w)}}  \right\}$$
$$\text{-- is first-order definable in $\mathbf{YF}^*$}. $$    
\end{Prop}
\begin{proof}
It is definable by the formula
    $$\phi_{d+}''(u,v,w):= \phi_{d+}'(u,v,w) \land \forall v'\forall v''\forall v'''(\phi_{d+}'(u,v',w)\land \phi_E(v'',v) \land \phi_E(v''',v')\to v''\geqslant v''' ) .$$
\end{proof}

\begin{Prop}    
$$\left\{(1^n,v,w)\in\mathbb{YF}^3: \quad d(w) \geqslant n\geqslant 0, w=1^{e_0}21^{e_1}2\ldots 2 1^{e_{d(w)-1}}21^{e_{d(w)}}, v=1^{e_n}\right\}$$
$$\text{-- is first-order definable in $\mathbf{YF}^*$}. $$    
\end{Prop}
\begin{proof}
It is definable by the formula
    $$\phi_{exp}(u,v,w):=  \exists v'\exists v''(\phi_{d+}''(u,v',w)\land \phi_E(v'',v') \land \phi_+(v'',v,u) ) .$$
\end{proof}

\begin{Prop}    
$$\left\{\left(u,1^{b(u)} \right)\in\mathbb{YF}^2 :\quad u=1^n, n\geqslant 2 \right\}\text{-- is first-order definable in $\mathbf{YF}^*$}. $$    
\end{Prop}
\begin{proof}
It is definable by the formula
    $$\phi_b'(u,v):= \phi_{1^n,11}(u)\land \phi_{1^n,11}(v) \land \forall v'\forall v'' \left(\psi_{Primexp}(v',v'',v) \to v'=\varepsilon \land r(u,v'') \right) .$$
\end{proof}

\begin{Prop}    
$$\left\{(1^n,v,w)\in\mathbb{YF}^3: \quad d(w) > n\geqslant 0, w=1^{e_0}21^{e_1}2\ldots 2 1^{e_{d(w)-1}}21^{e_{d(w)}}, v=1^{e_n}\right\}$$
$$\text{-- is first-order definable in $\mathbf{YF}^*$}. $$    
\end{Prop}
\begin{proof}
It is definable by the formula
    $$\phi_{exp}'(u,v,w):= \phi_{exp}(u,v,w) \land \forall w' (\phi_d(w',w)\to u\ngeqslant w') .$$
\end{proof}

\begin{Prop}    
$$\left\{(1^n,v,w)\in\mathbb{YF}^3: \quad d(w) = n\geqslant 0, w=1^{e_0}21^{e_1}2\ldots 2 1^{e_{d(w)-1}}21^{e_{d(w)}}, v=1^{e_{d(w)}+1}\right\}$$
$$\text{-- is first-order definable in $\mathbf{YF}^*$}. $$    
\end{Prop}
\begin{proof}
It is definable by the formula
    $$\phi_{exp}''(u,v,w):= \exists v'\exists w'( \phi_{exp}(u,v',w) \land  \phi_d(w',w)\land  u\geqslant w'  \land \phi_+(v,v',1)  ) .$$
\end{proof}

\begin{Prop}    
$$\left\{\left(u,1^{b(u)}\right)\in\mathbb{YF}^2: \quad d(u)\geqslant1\right\}\text{-- is first-order definable in $\mathbf{YF}^*$}. $$    
\end{Prop}
\begin{proof}
It is definable by the formula
   $$\phi_b''(u,v):= u\geqslant 2\land \forall v'\forall v'' \left(\psi_{Primexp}(v',v'',v) \longleftrightarrow \phi'_{exp}(v',v'',u) \lor  \phi''_{exp}(v',v'',u) \right) .$$
\end{proof}

\begin{Prop} \label{final}
$$\left\{\left(u,1^{b(u)}\right)\in\mathbb{YF}^2\right\}\text{-- is first-order definable in $\mathbf{YF}^*$}. $$    
\end{Prop}
\begin{proof}
It is definable by the formula
$$\phi_b(u,v):= id_\varepsilon(u)\land id_\varepsilon(v)\lor id_1(u)\land id_1(v)\lor \phi'_b(u,v) \lor \phi''_b(u,v).$$
\end{proof}

\begin{theorem}
$\left<\mathbb{YF},\geqslant, 2 \right>$ has the maximal definability property.

\end{theorem}

\newpage

\section{Acknowledgements}

This work was supported by the Ministry of Science and Higher Education of the Russian Federation (agreement 075-15-2025-344 dated 29/04/2025 for Saint Petersburg Leonhard Euler International Mathematical Institute at PDMI RAS).


\newpage

\addcontentsline{toc}{section}{References}


\begin{thebibliography}{}
    
\bibitem{Fo}
С. В. Фомин. \emph{Обобщенное
соответствие Робинсона -- Шенстеда -- Кнута}. Зап.
научн. сем. ЛОМИ, 155
(1986), 156-175.

\bibitem{St} R. P. Stanley. \emph{Differential posets}.
J. Amer. Math. Soc. 1 (1988), pp. 919-961.

\bibitem{Ma} Yuri Matiyasevich. \emph{Hilbert’s Tenth Problem}. MIT Press, Cambridge, 1993.

\bibitem{Tar} Alfred Tarski, Andrzej Mostowski, Raphael M. Robinson, \emph{Undecidable Theories}. North-Holland, Amsterdam, 1953.

\bibitem{W} Alexander Wires. \emph{Complexity in Young's lattice}.
Annals of Pure and Applied Logic 173 (2022), paper 103075.

\bibitem{L1} Simon Halfon, Philippe Schnoebelen, Georg Zetzsche. \emph{Decidability, complexity, and the expressiveness of first-order logic
over the subword ordering}, arXiv:1701.07470v1 [cs.LO], 2017.

\bibitem{L2} Jaroslav Ježek, Ralph Mckenzie. 
\emph{Definability in the lattice of equational theories of semigroups, I}, Semigroup Forum 46 (1993), pp. 199-245.

\bibitem{L3} Jaroslav Ježek, Ralph Mckenzie. \emph{Definability in substructure orderings I: finite semilattices}, Algebra Univers. 61 (2009), pp.
59-75.

 \bibitem{L4}Jaroslav Ježek, Ralph Mckenzie. \emph{Definability in substructure orderings II: finite ordered sets} Order 27 (2010), pp. 115-145.

\bibitem{L5} Jaroslav Ježek, Ralph Mckenzie. 
\emph{Definability in substructure orderings III: finite distributive lattices}, Algebra Univers. 61 (2009), pp. 283-300.

\bibitem{L6} Oleg V. Kudinov, Victor L. Selivanov, Lyudmilla V. Yartseva. 
\emph{Definability in the subword order}, in: Proc. CiE-2010, in:
LNCS, vol. 6158, Springer, 2010, pp. 246-255.

\bibitem{L7} Oleg V. Kudinov, Victor L. Selivanov. \emph{A Gandy theorem for abstract structures and applications to first-order definability},
in: Proc. of CiE-2009, in: LNCS, vol. 5635, Springer, Berlin, 2009, pp. 290-299.

\bibitem{L8} Oleg V. Kudinov, Victor L. Selivanov. \emph{Definability in the infix order on words}, in: Proc. of DLT-2009, in: LNCS, vol. 5583,
Springer, Berlin, 2009, pp. 454-465.

\bibitem{L9} Ádám Kunos. \emph{Definability in the embeddability ordering of finite directed graphs, II}, Order 36 (2019), pp. 291-311.

\end{thebibliography}
\end{document}